\documentclass[12pt]{article}
\usepackage{amsmath,amssymb,amsfonts,euscript}
\newtheorem{theorem}{Theorem}[section]
\newtheorem{lemma}[theorem]{Lemma}
\newtheorem{corollary}[theorem]{Corollary}
\newtheorem{proposition}[theorem]{Proposition}

\newtheorem{remark}[theorem]{Remark}
\title{Pure simplicial complexes and well-covered graphs}
\author{Rashid Zaare-Nahandi\\ Institute for Advanced Studies in Basic Sciences (IASBS),\\
Zanjan 45195, Iran \\ E-mail: rashidzn@iasbs.ac.ir}

\date{}

\begin{document}

\maketitle

\begin{abstract}
A graph $G$ is called well-covered if all maximal independent sets of vertices have the same cardinality.
A simplicial complex $\Delta$ is called pure if all of its facets have the same cardinality.
Let $\mathcal G$ be the class of graphs with some disjoint maximal
cliques covering all vertices. In this paper, we prove that for any simplicial complex or any graph, there is a corresponding graph in class $\mathcal G$ with the same well-coveredness property. Then some necessary and sufficient conditions are presented to recognize fast when a
graph in the class $\cal G$ is well-covered or not. To do this characterization, we use an algebraic interpretation according to zero-divisor elements of the edge rings of graphs.

\footnotetext{
\noindent{\it 2010 MR Subject Classification:} 05C25, 05E40, 05E45, 13F55.}
\end{abstract}
\section{Introduction}
A graph G is said to be well-covered (or unmixed) if every maximal independent sets of vertices
have the same cardinality. These graphs were introduced by M.~D.~Plummer~\cite{16} in
1970. Although the recognition problem of well-covered graphs in general is
Co-NP-complete~(\cite{19}), it is characterized for certain classes of graphs. For instance,
claw-free well-covered graphs~\cite{22}, well-covered graphs which have girth at least~5~\cite{7},
(4-cycle, 5-cycle)-free~\cite{8} or chordal graphs~\cite{18}
 are all
recognizable in polynomial time. Excellent surveys of works on
well-covered graphs are given in Plummer~\cite{17} and Hartnell~\cite{10}.

Let $G$ be a graph with no loop and multiple edge. Denote the set of vertices of $G$ by $V(G)$ and the set
of edges by $E(G)$. A subset $A$ of $V(G)$ is called an independent set if there is no any edge between vertices of $A$. Denote the cardinality of the largest independent set in $G$ by $\alpha(G)$.
A subset $C$ of $V(G)$ is called a clique if any two vertices in $C$ are adjacent.

Let $A$ and $B$ be subsets of $V(G)$. We say $A$ dominates $B$ if for any vertex $v$ in $B$, $v$ is in $A$ or there is at least one vertex in $A$ adjacent to $v$. The set $A$ is called a vertex cover of $G$ if any edge of $G$ has at least one edge in $A$. A vertex cover is called minimal if any proper subset of it is not a vertex cover.

A subset of $E(G)$ is called a matching if there is not any common vertex in any two edges in this set.
A matching is called perfect matching if it covers all vertices of $G$.

Let $[n] = \{1,2,\ldots,n\}$. A (finite) simplicial complex $\Delta$ on $n$ vertices, is a collection of subsets of
$[n]$ such that the following conditions hold:\\
a) $\{i\}\in\Delta$ for each $i\in [n]$,\\
b) if $E\in\Delta$ and $F\subseteq E$, then $F\in\Delta$.\\
An element of $\Delta$ is called a face and a maximal face with respect to inclusion is called a facet.
The set of all facets is denoted by $\mathcal{F}(\Delta)$.  The dimension of a
face $F\in\Delta$ is defined to be $|F|-1$ and dimension of $\Delta$ is maximum of dimensions of its faces.
A simplicial complex is called pure if all of its facets have the same dimension. For more details on simplicial complexes see \cite{Stan}.

Let $G$ be a graph. The set of all independent sets of vertices of $G$ is a simplecial complex, because,
any single vertex is independent and any subset of an independent
set is again independent. We assume that the empty set is also an independent set. This simplicial complex is called independence complex of $G$ and is denoted by $\Delta_G$.  With the above definitions,
a graph $G$ is well-covered means that the complex $\Delta_G$ is pure.

Let $\Delta$ be a simplicial complex on the vertex set $[n]$. The
barycentric subdivision of $\Delta$, denoted by $\mbox{bs}(\Delta)$, is
a simplicial complex with vertex set consisting of all nonempty
faces of $\Delta$. A face in $\mbox{bs}(\Delta)$ consists of comparable vertices, that is, two vertices
lie in a face in $\mbox{bs}(\Delta)$ if one is a subset of the other. In other words,
facets of $\mbox{bs}(\Delta)$ are maximal chains of faces of
$\Delta$ considered as a poset with respect to inclusion order.

It is easy to see that the minimal non-faces of $\mbox{bs}(\Delta)$ are
subsets of $\Delta$ with exactly two non-comparable elements.
Therefore, $\mbox{bs}(\Delta)$ is an independence  complex of a graph.
In fact this graph is non-comparability graph of $\Delta$. Vertices of the graph are nonempty faces of $\Delta$ and two vertices are adjacent if their corresponding faces are not comparable.
This graph is denoted by $G(\Delta)$.

It is
known that the dimension (and many other invariants) of a simplicial
complex and its barycentric subdivision are equal (\cite{BW} and \cite{KW}). Specially a simplicial complex $\Delta$ is pure if and only if its barycentric subdivision is pure and it is equal to say that the graph $G(\Delta)$ is well-covered.

\section{Well-covered graphs with clique covers}

Let $\cal G$ be the class of graphs such that for each $G\in \cal G$ there are $k=\alpha(G)$ cliques in $G$ covering all its vertices. Let $G\in \cal G$ and $Q_1,\ldots,Q_k$ be cliques such that $V(Q_1)\cup\cdots\cup V(Q_k)=V(G)$. In this case, we may take $Q'_1=Q_1$, and for $i=2,\ldots,k$, $Q'_i$ the induced graph on the vertices $V(Q_i)\setminus (V(Q_1)\cup\cdots\cup V(Q_{i-1}))$. Then, $Q'_1,\ldots,Q'_k$ are $k$ disjoint cliques covering all vertices of $G$. We call such a set of cliques, a basic clique cover of the graph $G$. Therefore, any graph in the class $\cal G$ has a basic clique cover. Note that, $k=\alpha(G)$ is the smallest number which the graph $G$ may has a clique cover. It is not true that any graph has a basic clique cover. For example, a cycle of length 4 has a basic clique cover consisting of 2 cliques but a cycle of length 5 does not have any basic clique cover.

\begin{proposition}\label{class}
Let $\Delta$ be a simplicial complex. Then, $G(\Delta)$ is in the class $\cal G$. Moreover, $\Delta$ is pure if and only if $G(\Delta)$ is well-covered.
\end{proposition}

{\it Proof}. Note that any two faces in $\Delta$ with the same dimension are not comparable. Therefore, for each $i$,  $0\leq i\leq \dim(\Delta)$, if $\Delta(i)$ is the set of all faces of $\Delta$ with dimension $i$, then there are not any two comparable face in this set and the corresponding vertices in the graph $G(\Delta)$ make a clique. These cliques are disjoint and cover all vertices of $G(\Delta)$. In fact the set of these cliques is a basic clique cover of $G(\Delta)$. The last statement is clear. \hfill $\Box$

Now, we give some criteria equivalent to well-covered property of graphs in the class $\cal G$.

\begin{theorem}\label{main} Let $G$ be a graph in the class $\cal G$ with a basic clique cover $Q_1, \ldots, Q_k$. Then $G$ is well-covered  if and only if for each $i$, $1\leq i\leq k$, if $A\subseteq V(G)\setminus Q_i$ dominates $Q_i$, then $A$ is not an independent set.
\end{theorem}

{\it Proof.} Assume that $G$ is well-covered. Let $1\leq i\leq k$ be given and $A\subseteq V(G)\setminus Q_i$ be a dominating set of $Q_i$. If $A$ is independent, then there is a maximal independent set $B$ containing $A$. But, $B\cap Q_i=\varnothing$ because any vertex of $Q_i$ is adjacent to some vertices in $A\subseteq B$. In other hand, $B$ has at most one element in common with each $Q_j$, $j\neq i$. Therefore, $|B|<k$ which is a contradiction with well-coveredness of $G$.

Conversely, let $A$ be a maximal independent set. Then $|A\cap Q_i|\leq 1$ for each $1\leq i\leq k$ and $|A|\leq k$.
The claim follows if one shows
$|A|=k$. So, if one assume $|A\cap Q_i|=\varnothing$ for some $i$, then one can apply the assumption
and $A$ is not dominating for $Q_i$, which means that there exists a $v \in Q_i$
not adjacent to any vertex of $A$. By maximality of $A$, $v\in A$ and hence $|A\cap Q_i|=1$
which is a contradiction to the assumption. So, finally one get $|A\cap Q_i|=1$ and
the claim follows. \hfill $\Box$

\begin{proposition}\label{partite} Let $G$ be a $s$-partite well-covered graph such that all maximal cliques are of size $s$. Then all parts have the same cardinality and there is a perfect matching between each two parts.
\end{proposition}

{\it Proof.} Let the $s$ parts of $G$ be $V_1, \ldots, V_s$. Let $1\leq i\leq s$ and $v\in V_i$. Each vertex  belongs to some maximal clique and each maximal clique intersects each part in exactly one vertex. Therefore, the vertex $v$ is adjacent to some vertices in each part $V_j$, $1\leq j\leq s$, $j\neq i$. Then the part $V_i$ is a maximal independent set because for each vertex out of $V_i$, there is an edge connecting it to some vertex in $V_i$. The graph $G$ is well-covered therefore, cardinality of parts are the same.

Let $1\leq i < j\leq s$ be two given integers.
Let $A\subseteq V_i$ be a nonempty set and $N_j(A)$ be the set of all vertices in $V_j$ adjacent to some vertices in $A$. Suppose $|N_j(A)|< |A|$. There is no any edge between $A$ and $V_j\setminus N_j(A)$. Therefore, $A\cup (V_j\setminus N_j(A))$ is an  independent set and its size is strictly greater than size of $V_j$, which is a contradiction with well-coveredness of $G$. Therefore, $|N_j(A)| \geq |A|$ for each nonempty subset $A$ of $V_i$. Therefore, by Theorem of Hall~\cite{Hall}, there is a set of distinct representatives (SDR) for the set $\{N_j(\{v\}) : v\in V_i\}$, which is a perfect matching between $V_i$ and $V_j$. \hfill $\Box$

\paragraph{Example.} It is not true that in a well-covered graph $G$, there are $\alpha(G)$ maximal cliques covering $G$. For instance, consider any cycle $C_n$ for odd $n$. In this case, $\alpha(C_n)=\frac{n-1}{2}$ and any $\frac{n-1}{2}$ cliques, which are edges, can not cover all vertices. Also the above statement is not true in class of all well-covered $s$-partite graphs. For instance consider the following graph which is 3-partite, well-covered with maximal independent sets of size 2. But, there are no two maximal cliques covering $V(G)$.

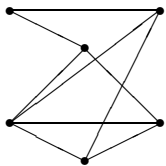
$$
\unitlength=1cm
    \begin{picture}(-4,-1)(3,1.5)
     \put(0,0.5){\circle*{0.1}}
     \put(1,0){\circle*{0.1}}
     \put(2,0.5){\circle*{0.1}}
     \put(0,2){\circle*{0.1}}
     \put(1,1.5){\circle*{0.1}}
     \put(2,2){\circle*{0.1}}
     \put(0,0.5){\line(2,-1){1}}
     \put(0,0.5){\line(1,0){2}}
     \put(0,0.5){\line(1,1){1}}
     \put(0,0.5){\line(4,3){2}}
     \put(1,0){\line(2,1){1}}
     \put(1,0){\line(1,2){1}}
     \put(2,0.5){\line(-1,1){1}}
     \put(0,2){\line(1,0){2}}
     \put(0,2){\line(2,-1){1}}
     %\put(1,2){\line(2,-1){1}}
     \end{picture}
$$

\vspace*{1.5cm}

Stating many examples motivates the following conjecture.

\paragraph{Conjecture.} Let $G$ be a $s$-partite well-covered graph with all maximal cliques of size $s$. Then, $G$ is in the class $\cal G$.

At the end of this section, we restate the result of Ravindra about well-covered bipartite graphs.

\begin{corollary}{\rm \cite{24}}  Let $G$ be a bipartite graph with no vertex of degree zero. Then, $G$ is well covered if and only if there is a perfect matching and for each $\{x,y\}$ in this matching, the induced subgraph on $N[\{x,y\}]$ is a complete bipartite graph.
\end{corollary}

{\it Proof.} Let $G$ be well-covered. By Proposition~\ref{partite}, cardinality of both parts are the same and there is a perfect matching in $G$. Moreovere, the edges in the matching make a basic clique cover of $G$. Let $\{x,y\}$ be an edge in the matching. By Theorem~\ref{main}, $G$ is well-covered if and only if any dominating set of $\{x,y\}$ is dependent. The last statement is equal to say that any vertex in $N(\{x\})$ is adjacent to any vertex in $N(\{y\})$,
i.~e., the induced subgraph on $N[\{x,y\}]$ is a complete bipartite graph. \hfill $\Box$

\section{An algebraic interpretation}
There is an interesting algebraic interpretation of well-coveredness of graphs in class $\cal G$, which we state in this section. First we recall some definitions in commutative algebra.

Let $G$ be a graph with vertex set $\{v_1, \ldots, v_n\}$ and $K$ be a field. In the polynomial ring $K[x_1, \ldots, x_n]$, consider $I(G)$ be the ideal generated by all monomials of the form $x_ix_j$ provided that $v_i$ and $v_j$ are adjacent in $G$. This ideal is called edge ideal of the graph $G$ and the quotient ring $R(G)=K[x_1,\ldots,x_n]/ I(G)$ is called edge ring of $G$. This ring is introduced by R. Villarreal \cite{Vil} and has been extensively studied by several mathematicians.

Let $R$ be a commutative ring. An element $a\neq 0$ in $R$ is called zero-divisor if there is a nonzero element $b\in R$ such that $ab=0$. An ideal in $R$ is called monomial ideal if it can be generated by a set of monomials. For example, edge ideal of a graph is a monomial ideal. In a ring of polynomials, it is well known and easy to check that a polynomial $f$ belongs to a monomial ideal if and only if each monomial of $f$ belongs to the ideal. If the monomial ideal is also square-free, then a monomial in $K[x_1,\ldots,x_n]$ belongs to $I$ if and only if its square-free part (its radical) belongs to $I$.
As an example of zero-divisor element, let $R(G)$ be the edge ring of a graph $G$. Let $v_i$ be adjacent to $v_j$ in $G$. The elements $x_i$ and $x_j$ are not zero in $R(G)$ but $x_ix_j=0$. Here, with abuse of notation, we have written $x_i$ as same as its image in $R(G)$.

A term ordering on $K[x_1,\ldots,x_n]$ is a linear order $\preceq$ on the set of terms $\{x_1^{a_1} x_2^{a_2} \ldots x_n^{a_n}\ : \ a_i\in {\Bbb Z}_{\geq 0}, i=1,2,\ldots n \}$, such that for each terms $\alpha, \alpha_1, \alpha_2$,
the following conditions hold.
\begin{itemize}
\item[a)]
if $\alpha_1 \preceq \alpha_2$ then $\alpha_1\alpha\preceq\alpha_2\alpha$.
\item[b)]
$1\preceq\alpha$.
\end{itemize}
Lexicographic, degree lexicographic and degree reverse lexicographic orders are examples of term orderings.  There is a rich literature about term orderings, for instance see \cite{Kr}.

\begin{lemma}\label{linearzero} Let $K$ be a field and $I\subseteq K[x_1,\ldots,x_n]$ be an ideal generated by square-free monomials. Let $f$ be a nonzero linear polynomial in $R=K[x_1,\ldots,x_n]/I$. Then, $f$ is zero-divisor in $R$ if and only if there is a nonzero square-free monomial  $m\in R$ such that $mf=0$.
\end{lemma}

{\it Proof.} Let $f$ be zero-divisor in $R$, then, there is a nonzero polynomial $g$ in $R$ such that $fg=0$.
We may rearrange variables such that $f=x_1+a_2x_2+\cdots+a_sx_s$, $a_j\in K$. Let $\prec$ be the lexicographic order on terms of $K[x_1,\ldots,x_n]$ with respect to $x_1\succ x_2 \succ \cdots \succ x_n$. Let  $g=m_1 + m_2 + \cdots + m_t$ be decomposition of $g$ to nonzero monomials such that $m_1 \succ m_2 \succ\cdots \succ m_t$. Then, in $fg$, the monomial $x_1m_1$ is strictly greater than all other monomials. Therefore, $x_1m_1$ must be zero in $R$. The ideal $I$ is square-free and $x_1m_1\in I$, therefore, we may assume that $x_1\nmid m_1$. By the lexicographic order, we have $x_1\nmid m_i$ for all $1\leq i\leq t$. In other hand, $fg-x_1m_1\in I$. The greatest term of $fg-x_1m_1$ is $x_1m_2$ and then $x_1m_2\in I$ and $fg-(x_1m_1+x_1m_2)\in I$. Continuing this process, we have $x_1m_i\in I$ for all $1\leq i\leq t$ and therefore, $fg-x_1g\in I$. In the polynomial $fg-x_1g$ the greatest term is $x_2m_1$ which must be in $I$. Similarly, $x_2m_i\in I$ for all $1\leq i\leq t$. Finally, we get $x_im_j\in I$ for each $1\leq i\leq s$ and  $1\leq j\leq t$. It means that $m_if\in I$ for each $1\leq i\leq t$. Specially $m_1f\in I$, and because $I$ is square-free and $f$ is linear, we may take  $m_1$ to be square-free. The converse is trivial by definition. \hfill $\Box$\\

Note that in the above lemma, assuming that $I$ is square-free is essential. Because, for example in $K[x_1,x_2]$ assume that $I=\langle x_1^3, x_2^3\rangle$. Then, $(x_1-x_2)(x_1^2+x_1x_2+x_2^2)\in I$, that is $(x_1-x_2)$ is a zero-divisor in  $K[x_1,x_2]/I$ but, there is no any nonzero square-free monomial eliminating $(x_1-x_2)$.

\begin{theorem}\label{zerodiv} Let $G$ be a graph in the class $\cal G$ and $\alpha(G)=k$. Let  $Q_1, \ldots, Q_k$ be a basic clique cover of $G$. Consider
$$
\theta_i=\sum_{v_j \in Q_i} x_j, \ \ \ \ \ i=1,\ldots,k.
$$
Then, $G$ is well-covered if and only if for each $i=1,\ldots,k$, the polynomial $\theta_i$ is not zero-divisor in the ring $R(G)$.
\end{theorem}

{\it Proof.} Let $\theta_i$ be zero-divisor in $R(G)$. By Lemma~\ref{linearzero}, the polynomial $\theta_i$ is zero-divisor in $R(G)$ if and only if there is a nonzero square-free monomial $m$ in $R(G)$ such that $m\theta_i=0$ or equivalently $m\theta_i\in I(G)$.  The ideal $I(G)$ is a monomial ideal, then, for each $v_j$ in $Q_i$, we have $mx_j\in I(G)$. Let $m=x_{i_1}\cdots x_{i_r}$ and $A=\{v_{i_1}, \ldots, v_{i_r}\}$. Then,
$mx_j\in I(G)$ means that there is a vertex $v_{i_l}$ in $A$ such that $v_{i_l}$ is adjacent to $v_j$. This means that the set $A$ is a dominating set of $Q_i$. In other hand, if $v_j$ is in $A\cap Q_i$, then $x_j\theta = x_j^2$ in $R(G)$ and there is $v_{i_l}$ in $A$ adjacent to $v_j$ and therefore $m=0$ in $R(G)$ which is a contradiction. Therefore $A\subseteq V(G)\setminus V(Q_i)$.
Note that $A$ is independent if and only of $m$ is not zero in $R(G)$. Now, Theorem~\ref{main} implies that if $\theta_i$ is a zero-divisor in $R(G)$ for some $1\leq i\leq k$, then, $G$ is not well-covered.

Conversely, if $G$ is not well-covered then, again by Theorem~\ref{main}, there is an independent set $\{v_{i_1}, \ldots, v_{i_r}\}\subseteq V(G)\setminus V(Q_i)$ which dominates $Q_i$ for some $1\leq i\leq k$. In this case, $m=x_{i_1}\cdots x_{i_r}$ is a nonzero monomial in $R(G)$ such that $m\theta_i=0$ and $\theta_i$ is zero-divisor. This completes the proof. \hfill $\Box$\\

Let $G$ be a graph in the class $\cal G$. Then, by Theorem~\ref{zerodiv}, $G$ is well-covered if and only if each polynomial $\theta_i$ is non-zero-divisor in the ring $R(G)$. In other hand, the set of all zero-divisors of $R(G)$ is union of all minimal primes of the ideal $I(G)$. Minimal primes of $I(G)$ are corresponding to minimal vertex covers of $G$. Therefore, checking well-coveredness of the graph $G$ is equal to check that for each $i$, $1\leq i\leq k$, the set of vertices of $Q_i$ is a part of a minimal vertex cover of $G$ or not. But, this is a simple task: it is enough to check that the set of vertices of $Q_i$ is a minimal vertex cover of the induced sub-graph of $G$ on $N(Q_i)$, which can be done in a polynomial time algorithm. Therefore, we have proved the following. 

\begin{corollary}
The well-coveredness of a graph in the class $\cal G$ can be checked in polynomial time. 
\end{corollary}

We know that an arbitrary graph $G$ is well-covered if and only if the corresponding graph $G(\Delta_G)$ is well-covered.   The graph $G(\Delta_G)$ is in the class $\cal G$ and its well-coveredness can be checked in polynomial time. But, this does not solve completely the problem of well-covered checking of graphs, because passing from $G$ to $G(\Delta_G)$ can not be done in polynomial time. In fact,  the graph $G(\Delta_G)$ has a huge number of vertices in comparison with $G$.

\begin{remark} \rm The next natural question is when a graph in the class $\cal G$ is Cohen-Macaulay. With the notations above, Cohen-Macaulayness of $G$ is equal to regularity of the sequence $\theta_1, \theta_2, \ldots, \theta_k$ in $R(G)$. It means that $\theta_1$ is not zero-divisor in $R(G)$ and for $i=2, \ldots, k$, the element $\theta_i$ is not zero-divisor in $R(G)/\langle \theta_1, \ldots, \theta_{i-1}\rangle$. Therefore, one can say that if $G$ is Cohen-Macaulay, then $G\setminus Q_i$ is Cohen-Macaulay for each $1\leq i\leq k$. It is well-known that a simplicial complex $\Delta$ is Cohen-Macaulay if and only if the graph $G(\Delta)$ is Cohen-Macaulay~(\cite{Stan}). Therefore, to check Cohen-Macaulayness of all simplicial complexes and all graphs, it is enough to check Cohen-Macaulayness of all graphs in the class $\cal G$.
\end{remark}

\end{document}